\documentclass[12pt]{article}

\usepackage{amssymb}
\usepackage{amsmath}
\usepackage{amsthm}
\usepackage{hyperref}

\newtheorem{theorem}{Theorem}[section]

\newtheorem{proposition}[theorem]{Proposition}

\theoremstyle{definition}

\theoremstyle{remark}

\numberwithin{equation}{section}

\newcommand{\vp}{\varphi}

\newcommand{\la}{\langle}
\newcommand{\ra}{\rangle}

\newcommand{\tos}{\rightrightarrows}

\newcommand{\rn}{\mathbb{R}}

\newcommand{\HH}{\mathcal{H}}  
  
\newcommand{\JJ}{\mathcal{J}}

\title{Maximal monotonicity, conjugation and the duality product\\
{\small published on:
  \href{http://www.ams.org/proc/2003-131-08/S0002-9939-03-07053-9/home.html}{
         {\it Proc.\ Amer.\ Math.\ Soc.\ } {\bf 131} (2003) 2379--2383.
         }
}
}
\author{ {\it Regina Sandra Burachik}\thanks{Partially supported
                     by CNPq
               and   by PRONEX--Optimization.}
         \\Engenharia de Sistemas e Computa\c c\~ao
         \\  COPPE--UFRJ CP 68511
         \\  Rio de Janeiro--RJ
         \\  CEP 21945--970 Brazil
          regi@cos.ufrj.br
\and
         {\it B. F. Svaiter}\thanks{Partially supported by 
                    CNPq Grant 301200/93-9(RN)
                    and by PRONEX--Optimization.}
         \\ IMPA Instituto de Matem\' atica Pura e Aplicada
         \\ Estrada Dona Castorina 110
         \\ Rio de Janeiro--RJ
         \\ CEP 22460-320 Brazil
          benar@impa.br
}

\date{17 June 2002}

\begin{document}

\maketitle

\begin{abstract}
Recently, the authors studied the connection between each maximal
monotone operator $T$ and a family $\HH(T)$ of convex functions.
Each member of this
 family characterizes the operator and  satisfies two particular
inequalities.

The aim of this paper is to establish the converse of the latter fact.
Namely, that every convex function satisfying those two particular
inequalities is associated to a unique maximal monotone operator.
\\
2000 Mathematics Subject Classification:  47H05
\\
{ \sc keywords:} convex functions, maximal monotone operators, duality product, 
          conjugation
\end{abstract}


\section{Introduction}
Let $X$ be a real Banach space and $X^*$ be the dual of $X$.
Denote by $\la\cdot,\cdot\ra$ the duality product.
A multivalued operator $T:X\tos X^*$ is \emph{monotone} if
\[ \la x_1-x_2,v_1-v_2\ra\geq 0, \quad \forall v_1\in T(x_1), v_2\in T(x_2).\]
Such an operator is \emph{maximal monotone} if its \emph{graph},
that is, the set
\[ G(T)=\{ (x,v)\in X\times X^*\;|\; v\in T(x)\}, \]
is not properly contained in the graph of any other monotone operator
$T':X\tos X^*$. We will identify $T$ with its graph $G(T)$.

The \emph{subdifferential} of a function $f:X\to \overline{\rn}$ is the 
multivalued operator $\partial f:X\tos X^*$ defined by
\[ \partial f(x)=\{v\in X^*\,|\, 
  f(y)\geq f(x)+\la y-x,v\ra,\,\forall y\in X\}.\]
A convex function is \emph{closed} if it
is lower semicontinuous and it is \emph{proper} if it not attains
the value  $-\infty$, and is not $+\infty$ everywhere.

Rockafellar \cite{rockafellar:1970mmsub} proved that subdifferentials
of proper closed convex functions on $X$ are maximal monotone.  In
general, maximal monotone operators are not subdifferentials
of convex functions. 
Krauss~\cite{krauss} managed to represent maximal monotone operators
by subdifferentials of saddle functions on $X\times X$.
After that, Fitzpatrick~\cite{fitzpatrick:1988} proved that maximal monotone
operators can be represented by convex functions on $X\times X^*$.
This result has been recently rediscovered in
\cite{martinez-thera:2001,burachik-svaiter:2002}.

Next we describe
Fitzpatrick's results.  Given a maximal monotone operator $T:X\tos
X^*$, define
\begin{equation}
  \label{eq:defh}
 \HH(T):=
\left\{
\begin{array}{l|l}
 & h\mbox{ convex, closed};\\
h:X\times X^*\to \overline{\rn} 
 & \forall (x,v)\in X\times X^*, h(x,v)\geq \la x,v\ra,\;\\
 & (x,v)\in T\Rightarrow h(x,v)=\la x,v\ra.
\end{array}
\right\}\;. 
\end{equation}
Define also $\vp_T:X\times X^*\to \overline{\rn}$,
\[ \vp_T(x,v):=\sup\{ \la x-y,u-v\ra\;|\; (y,u)\in T\}+\la x,v\ra.\]
\begin{theorem}[Fitzpatrick\cite{fitzpatrick:1988}]\label{th:f1}
Let $T:X\tos X^*$ be maximal monotone.
The function $\varphi_T$ belongs to $\HH(T)$ and is the smallest function
of
this family.
Moreover, for any $h\in\HH(T)$, 
\[ (x,v)\in T\Leftrightarrow h(x,v)=\la x,v \ra.\]
\end{theorem}
\noindent
From the above equivalence, it follows that each $h\in\HH(T)$ fully
characterizes $T$.

Given a function $f:X\to\overline{\rn}$, the \emph{Legendre transform}
or \emph{conjugate} of $f$ is defined as $f^*:X^*\to\overline{\rn}$,
\[ f^*(v)=\sup\{ \la x,v\ra-f(x)\;|\; x\in X\}.\]
Conjugation is an essential tool in the study of convex functions.
Let $f$ be a proper convex function. 
From the previous definitions, we have the \emph{Fenchel--Young inequality:}
for all $x\in X$, $v\in X^*$
\[f(x)+f^*(v)\geq \la x, v\ra\, ,\:\:
f(x)+f^*(v)= \la x, v\ra \iff v\in\partial f(x).
\]

For $h:X\times X^*\to\overline{\rn}$, the conjugate of $h$, is defined
on $X^*\times X^{**}$.  Since there is a natural injection of $X$ in to
$X^{**}$, we define $\JJ(h):X\times X^*\to\overline{\rn}$,
\[\begin{array}{rcl} 
  \JJ(h)(x,v) & = & h^*(v,x), \\
     & = & \sup\left\{ 
                  \left\la \left(y,u\right),\left(v,x\right)
                  \right\ra
                  -h(y,u);|\; (y,u)\in X\times X^*
                \right\},\\[.2em]
     & = & \sup\left\{ 
                  \la y,v\ra +\la x,u\ra-h(y,u);|\; (y,u)\in X\times X^*
               \right\}.
   \end{array}
\]
Fitzpatrick \cite[Prop.\ 4.2]{fitzpatrick:1988} proved that if $T:X\tos X^*$ is maximal
monotone, $\JJ(\vp_T)$ also belongs to $\HH(T)$. 
In \cite{burachik-svaiter:2002} this result was extended to any $h\in\HH(T)$.
Namely, if 
$h\in\HH(T)$, then $\JJ(h)\in\HH(T)$. 
Altogether, the result in \cite{burachik-svaiter:2002} can be expressed
as the implication 
\[
\begin{array}{l}
T:X\tos X^*\mbox{ maximal monotone}\\
h\in\HH(T)
\end{array}
\Rightarrow
\begin{array}{l}
\forall (x,v)\in X\times X^*,\\
h(x,v)\geq \la x,v\ra,\;h^*(v,x)\geq \la x,v\ra.
\end{array}
\]
Our aim is to prove the converse of this implication in a \emph{reflexive}
Banach space. Namely, 
\[
\begin{array}{l}
h:X\times X^*\to \overline{\rn} \hbox{ convex, lsc},\\ 
\forall (x,v)\in X\times X^*,\\
\;h(x,v)\geq \la x,v\ra,\;h^*(v,x)\geq \la x,v\ra.
\end{array}
\Rightarrow
\begin{array}{l}
\exists!\,T:X\tos X^*\mbox{ maximal monotone}\\
h\in\HH(T)
\end{array}
\]

The paper is organized as follows. In Section \ref{sec:te}
we state some necessary previous results.
The last section contains the formal statement and the proof of 
the implication above (see Section \ref{sec:main}, Theorem \ref{th:main}).

\section{Theoretical Preliminaries}\label{sec:te}

We include in this section theoretical results which
are necessary for the proof of Theorem \ref{th:main}. 
From now on  $X$ is a real Banach space. 
\begin{theorem}[{\cite[Theor. 5.3]{burachik-svaiter:2002}}]\label{th:jh}
Let $T:X\tos X^*$ be maximal monotone. Then, the operator $\JJ$ maps 
$\HH(T)$ into itself.  
\end{theorem}

We also assume from now on that $X$ is  \emph{reflexive}. 
Asplund \cite{asplund:1967} has shown that, in this case, there exists an
equivalent norm on $X$ which is everywhere G\^ ateaux differentiable except
at the origin and whose polar norm on $X^*$ is
everywhere G\^ ateaux differentiable except at the origin.  For
simplifying the notation, we assume that the given norm on $X$ already
has these special properties. We use the same notation $\|\cdot\|$ 
for this norm on $X$ and its associated norm on the dual $X^*$. 
Denote by $J$ the G\^ ateaux
gradient of the function $g(x)=(1/2)\|x\|^2$.  Thus, $J$ is the
duality mapping, which assigns to each $x\in X$ the unique $J(x)\in
X^*$ such that
\begin{equation}
  \label{eq:rock.2.1}
  \la x,J(x)\ra=\|x\|^2=\|J(x)\|^2.
\end{equation}
The inverse of this duality mapping will be denoted by $J_*$, which is
the subgradient of the function $g^*(v)=(1/2)\|v\|^2$.

The following result was proved in \cite[Section 2]{rockafellar:1970ms}, 
where it
appears as a corollary.
\begin{proposition}\label{pr:rock}
 Let
$T:X\tos X^*$ be a monotone operator.
Under the above assumptions,
 in order that $T$ be maximal monotone,
it is necessary and sufficient that $(T+J):X\tos X^*$ be onto. 
\end{proposition}

Using the fact that $\la z,u\ra\geq -\|z\|\,\|u\|$ for all $z\in X,\,u\in X^*$,
one can easily obtain the proposition below. 
\begin{proposition}\label{pr:el}
Under the above assumptions, let $z\in X$, $u\in X^*$. Then 
\[ \|z\|^2+\|u\|^2+2\la z,u\ra\geq 0, \]
with equality if and only if
\( u=-J(z)\) (or equivalently
\(z=-J_*(u)\) ).
\end{proposition}

\section{Main Result}\label{sec:main}
Now we state formally and prove the main result.  Recall that $X$ is a
reflexive real Banach space.  For convenience, $X$ has been Asplund-renormed.

\begin{theorem}\label{th:main}
Under the above assumptions, let 
$h:X\times X^*\to \overline{\rn}$ be a 
convex lower semicontinuous function.
Suppose that
\begin{equation}
  \label{eq:th.main.hyp}
  \forall (x,v)\in X\times X^*,\quad
  h(x,v)\geq \la x,v\ra,\;\; h^*(v,x)\geq \la x,v\ra.
\end{equation}
Define
\begin{equation}
  \label{eq:th.main.deft}
  T=\{ (x,v)\in X\times X^*\;|\; h(x,v)=\la x,v \ra\}.
\end{equation}
Then $T$ is maximal monotone
and $h$, $\JJ(h)\in\HH(T)$.
\end{theorem}

\begin{proof}
First we claim that
 $T$ is monotone. Indeed, take $v_1\in T(x_1)$, 
$v_2\in T(x_2)$, then
\begin{equation}
  \label{eq:tc1}
  \la x_1,v_1\ra=h(x_1,v_1),\quad \la x_2,v_2\ra=h(x_2,v_2).
\end{equation}
The convexity of $h$ together with \eqref{eq:th.main.hyp} gives
\begin{equation}
  \label{eq:tc2}
  \begin{array}{rcl}
(1/2)(h(x_1,v_1)+h(x_2,v_2))&\geq& h( (1/2)(x_1+x_2),(1/2)(v_1+v_2))\\
 &\geq& \la (1/2)(x_1+x_2),(1/2)(v_1+v_2)\ra
\end{array}
\end{equation}
Combining this with \eqref{eq:tc1} we obtain
\[ (1/2)\bigg( \la x_1,v_1\ra+\la x_2,v_2\ra\bigg)
\geq (1/4)\la x_1+x_2,v_1+v_2\ra,
\]
which is equivalent to $\la x_1-x_2,v_1-v_2\ra\geq 0$.

Now we claim that $T+J:X\tos X^*$ is onto.  To prove this fact, take
an arbitrary $v_0\in X^*$ and define
$\varphi:X\times X^*\to \rn\cup\{+\infty\}$,
\begin{equation}
  \label{eq:defphi}
\begin{array}{rcl}
   \varphi(x,v)&:=&(1/2)\bigg(\|x\|^2+\|v-v_0\|^2+2\la x,v-v_0\ra\bigg) +
                       \bigg(h(x,v)-\la x,v\ra\bigg)\\[.5em]
  &=&(1/2)\bigg(\|v-v_0\|^2+\|x\|^2\bigg)-\la v_0 ,x\ra +h(x,v),
\end{array}
\end{equation}
where $\|\cdot\|$, $J$ are the norm and duality map defined above,
respectively.  By the first expression for $\varphi$, assumptions
\eqref{eq:th.main.hyp}-\eqref{eq:th.main.deft} and Proposition
\ref{pr:el}, we have $\vp\geq 0$, with equality only if $v-v_0=-J(x)$
and $v\in T(x)$. This
implies $v_0\in (T+J)(x)$.
The second expression of $\vp$ shows that this function is the sum of
a differentiable convex function plus a lower semicontinuous convex
function.  By \cite[Th. 3]{rockafellar:1966} or
\cite[p. 62]{moreau:1967} the subdifferential of this sum is the sum of
the subdifferentials. Using also the equalities 
$J(\cdot)=\partial(1/2\|\cdot \|^2)$ and 
$J_*(\cdot)=\partial(1/2\|\cdot \|^2)$ we obtain
\[
\begin{array}{rcl}
\partial\vp(x,v)&=&(\partial_X \vp(x,v), \partial_{X^*} \vp(x,v)) + 
\partial h(x,v)\\
&=&(J(x)-v_0,J_*(v-v_0))+\partial h(x,v).
\end{array}
\]
Since $X\times X^*$ is reflexive and $\varphi$ 
is lower semicontinuous and strongly convex, it
attains a minimum at some $(x,v)
\in X\times X^*$.  Hence, for such $(x,v)$
\[  0\in (J(x)-v_0,J_*(v-v_0))+\partial h(x,v).\]
To simplify the manipulations, define
\begin{equation}
  \label{eq:drr}
  \begin{array}{rcl}
r&=&J_*(v-v_0)+x,\\
\rho&=&v-v_0+J(x).
  \end{array}
\end{equation}
With this notation, the last inclusion becomes
\[ (v-\rho,x-r)\in\partial h(x,v).\]
Hence, by Fenchel-Young we have that
\begin{equation}
  \label{eq:fy}
  h(x,v)+h^*(v-\rho,x-r)=\la x,v-\rho\ra+\la x-r,v\ra.
\end{equation}
Define now
\begin{equation}\label{eq:dc}
C:=\bigg( \la    x,   v\ra -h(    x,   v)\bigg) +
  \bigg( \la   v-\rho, x-r \ra -h^*(   v-\rho,x-r)\bigg).
\end{equation}
Assumption \eqref{eq:th.main.hyp} yields $C\leq 0$.
Using now \eqref{eq:fy},\eqref{eq:drr} we obtain
\[ 
\begin{array}{{rcl}}
C&=& \bigg( \la    x,   v\ra+\la   v-\rho, x-r \ra\bigg)
  -\bigg( \la x,v-\rho\ra+\la x-r,v\ra\bigg)\\
 &=&\la r,\rho\ra\\
 &=&\la x,J(x)\ra+\la J_*(v-v_0),v-v_0\ra +\la x,v-v_0\ra
                        +\la J_*(v-v_0),J(x)\ra.
\end{array}
\]
Using \eqref{eq:rock.2.1} and the fact that $J_*$ is the inverse of $J$ we have
\[\begin{array}{lcl}
 \la x,J(x)\ra & = & (1/2)\bigg(\|x\|^2+\|J(x)\|^2\bigg),\\[.4em]
\la J_*(v-v_0),v-v_0\ra & = & (1/2)\bigg(\|v-v_0\|^2+\|J_*(v-v_0)\|^2\bigg).
\end{array}
\]
The combination of 
these equalities with the above expression for $C$ gives
\begin{equation}\label{eq:c2}
\begin{array}{rcl}
 C&=&(1/2)\bigg( \|x\|^2+\|v-v_0\|^2+2\la x,v-v_0\ra\bigg)\\
  &&{}
+(1/2)\bigg(\|J(x)\|^2+\|J_*(v-v_0)\|^2+2\la J_*(v-v_0),J(x)\ra\bigg).
\end{array}
\end{equation}
By the above equation and  Proposition \ref{pr:el}, $C\geq 0$.
Therefore, $C=0$. Using this fact and again \eqref{eq:c2} and
Proposition \ref{pr:el}, we conclude that $v-v_0=-J(x)$, that is
\begin{equation}
  \label{eq:eq}
  v+J(x)=v_0.
\end{equation}
On the other hand, using \eqref{eq:dc}, \eqref{eq:th.main.hyp} and the
equality $C=0$ we conclude that $h(x,v)=\la v,x \ra$, which yields
\[ v\in T(x).\]
Therefore, $v_0\in (T+J)(x)$. Since $v_0$ is arbitrary, $T+J$ is onto.

We have thus proved that $T$ is monotone and $T+J$ is onto, hence by
Proposition \ref{pr:rock}, $T$ is maximal monotone.  It remains to
prove that $h$ and $\JJ(h)\in\HH(T)$.  In order to do this, we use
\eqref{eq:th.main.hyp} and the definition of $\HH$, for concluding
that $h\in\HH(T)$. The inclusion $\JJ(h)\in\HH(T)$ now follows from
Theorem \ref{th:jh}.
\end{proof}

\end{document}